\input amstex

\def\bb#1#2{b^{#1}_{#2}}
 \def\aaa#1#2{a^{#1}_{#2}}
\def\bcc#1{\Bbb C^{#1}}
\def\brr#1{\Bbb R^{#1}}

\def\bccc{\Bbb C }\def\brrr{\Bbb R }

\def\hd{, \hdots ,}
\def\inv{{}^{-1}}

\def\pp#1{\Bbb P^{#1}}
\def\ppp{\Bbb P }

\def\upp#1{p_{#1}}
\def\lp#1{p_{#1}}
\def\uup#1#2{p^{#1}_{#2}}

\def\qq#1{q_{#1}}
\def\qqq#1#2{q^{#1}_{#2}}

\def\tdim{\text{dim}}
\def\tcodim{\text{codim}}

\def\tspan{\text{Span}}
\def\ttrace{\text{trace}}

\def\upt{{}^t\!}

\def\xx#1{x^{#1}}

\def\zz#1{z_{#1}}

\documentstyle{amsppt}
\magnification = 1200
\hsize =15truecm
\hcorrection{.5truein}
\baselineskip =18truept
\vsize =22truecm
\NoBlackBoxes
\topmatter
\title
On an unusual
conjecture of Kontsevich and  variants of
Castelnuovo's lemma \endtitle
\rightheadtext{Kontsevich's unusual conjecture}
\author
  J.M. Landsberg
\endauthor

\address{Department of Mathematics,
Columbia University, New York,  NY 10027}
\endaddress
\email {jml\@math.columbia.edu}
\endemail
\date {   July 30, 1996 }\enddate
\thanks {Supported   by NSF grant
DMS-9303704.}
\endthanks
\endkeywords
\abstract{Let $A=(a^i_j)$ be an orthogonal matrix
(over $\brrr$ or $\bccc$)
with no entries zero.  Let $B=(b^i_j)$ be the matrix defined by
$b^i_j=\frac 1{a^i_j}$. M. Kontsevich conjectured that the
rank of $B$ is never equal to three. We interpret this
conjecture   geometrically and
prove it. The geometric statment can be understood as
 variants of the Castelnouvo lemma
and Brianchon's theorem.}
\endabstract

\endtopmatter
\document

\noindent\S 1. {\it Definitions and Statements}

\smallpagebreak

\noindent{\bf Definition 1.1}.
Given a $k\times l$ matrix $A= (a^i_{\alpha})$,
$1\leq i\leq k, 1\leq \alpha\leq l$,
with no entries zero, define
the {\it Hadamard inverse of $A$}, $B= (b^i_{\alpha})$, by $b^i_{\alpha}
=\frac 1{a^i_{\alpha}}$. (The name is in analogy with the
Hadamard product.)

Maxim Kontsevich conjectured the following:

\proclaim{Conjecture 1.2}(Kontsevich (1988))
Let $A$ be an orthogonal matrix
(over $\brrr$ or $\bccc$)
with no entries zero.
 Let $B$ be the Hadamard inverse of $A$. Then the
rank of $B$ is never equal to three.
\endproclaim

 At first glance, (1.2) may not appear all that striking because
based on a naive count,
one would not expect any low rank
Hadamard inverses of orthogonal matrices to exist
(see (1.16)).
However Kontsevich asserted and we   show the following:

\proclaim{Theorem 1.3} The space of $m\times m$ orthogonal matrices with
rank two Hadamard inverses is $(2m-3)$-dimensional.
\endproclaim

 We will rephrase (1.2), (1.3)   in geometric language and
prove them. First we will need some definitions:

\smallpagebreak

\noindent{\bf Definitions 1.4}.
Let $V=\bcc{n+1}$ or $\brr{n+1}$, let $  Q\in S^2V^*$ be a nondegenerate
quadratic form. Two points
  $z,w\in \ppp V$ are said to be {\it polar}, or more
precisely, {\it Q-polar}
if the corresponding
lines in $V$ are $  Q$-orthogonal. Given a point $z\in \ppp V$, the
{\it polar hyperplane of $z$ with respect to $Q$}, $H_{Q,z}$,
is   the hyperplane  of $Q$-polar
points to $z$.
A set of points $\zz 0\hd \zz n$, none lying
on the quadric in $\ppp V$ defined by $Q$, is said to be {\it apolar}
with respect to $Q$ if they are all mutually polar. (In particular,
vectors in $V$ representing them give a $  Q$-orthogonal basis of $V$.)
A quadratic form $  P$ is said to be apolar to $  Q$ if
 trace${}_{  Q}P=0$.

\smallpagebreak

 In what follows,
unless we specify the ground field, it may be taken to be
$\brrr$ or $\bccc$.

\smallpagebreak

\noindent{\bf Definition 1.5}.
Given a set of points $\zz 0\hd \zz n$ spanning $\pp n$,
define the {\it Cremona transform of $\pp n$
with respect to $\{ \zz i \}$} to be the
rational map
 $\phi :\pp n\dasharrow \pp n$,
 obtained by first blowing up the
codimension two spaces spanned by subsets of
$(n-1)$-tuples of the points,
 and then blowing down the $n$ hyperplanes containing
sets of $n$ of the points.
In coordinates, if
$$
\zz i = [0\hd 0,1,0\hd 0 ]\tag 1.6
$$ where the
$1$ occurs in the $i$-th slot,   the map is
$$
\phi ( [\xx 0\hd \xx n]) = [\frac 1{\xx 0}\hd \frac 1{\xx n}].\tag 1.7
$$
Note that the images of the blown down hyperplanes determine
a coordinate simplex in the image $\pp n$ which we will call
the {\it image simplex}.

Conjecture (1.2) is equivalent to:
\proclaim{Theorem 1.8, Version 1} Let $\zz 0\hd\zz n$
and $\upp 0\hd\upp n$ be two sets
 of points
spanning $\pp n$, each set  apolar
with respect to a nondegenerate quadratic form $Q$.
Let $\phi$ denote the Cremona transform defined by the
$\{\zz i\}$. If   the points $\phi (\upp i)$
fail to span a $\pp 3$, then in fact they span exactly a
$\pp 1$.

\endproclaim

\proclaim{Theorem 1.8, Version 2}
 Let $\zz 0\hd\zz n$
and $\upp 0\hd\upp n$ be two sets
 of points
spanning $\pp n $, each set  apolar
with respect to a nondegenerate quadratic form $Q$.
 Let $\Sigma $ be the
space of hypersurfaces of degree $n$ having multiplicities
of order $n-2$ at each $\zz i$. ($\Sigma$ is a $\pp{n}$.)
Let $\Lambda:=\{ P\in  \Sigma\, |\,  \upp i\in P\ \forall i \}$.
If $\tcodim \Lambda \leq 2$ then $\tcodim \Lambda =1$.
In this case the points $\zz i ,p_i$ all lie on a rational
normal curve.
\endproclaim

\demo{Equivalence of versions 1 and 2} $\Sigma$ is   the
space of
inverse images under $\phi$ of the   hyperplanes in
$\pp n$, and for any $\pp 1\subset \pp n$, $\phi\inv (\pp 1)$ is
a (possibly degenerate) rational normal curve. (The $\pp 1$'s that
we will consider will yield non-degenerate rational normal
curves.)\enddemo

One form of the classical Brianchon theorem says that given a conic
in the plane and given two triangles circumscribing the conic,
then the six points consisting of the vertices of the two triangles
all lie on another conic. The $n=3$
case of (1.8) provides  variant of this
which  was originally proven by Weddle
and Zeuthen (see [EP2] for an exposition of their work).

 \proclaim{Corollary 1.9 (Weddle, Zeuthen)} Let $Q\subset\pp 2$ be a smooth
conic.
Let $\zz 0,\zz 1,\zz 2$ and $\lp  0,\lp 1, \lp 2$ be two sets
of  $Q$-apolar points. Then the six
points $\zz 0,\zz 1,\zz 2$, $\lp  0,\lp 1, \lp 2$ all lie on a conic.
\endproclaim

Writing $\upp j=[\uup 0j,\uup 1j,\uup 2j]$,
  $z_i$ as in (1.6), and let $Q$ have equation $ \Sigma_i(\xx i)^2=0$,
 the conic all six points lie on is
$$
(\uup 20\uup 21\uup 22)\xx 0\xx 1 +
(\uup 10\uup 11\uup 12)\xx 0\xx 2 +
(\uup 00\uup 01\uup 02)\xx 1\xx 2. \tag 1.10
$$

More generally, we have:

\proclaim{Proposition 1.11}
 Let $\zz 0\hd\zz n$
and $\upp 0\hd\upp n$ be two sets
 of points
spanning $\bccc\pp n=\ppp V$, each set  apolar
with respect to a smooth quadric $Q$.
 Let
$\Cal Q = \{ P\in \ppp S^2V^* \,|\,  \zz i,\upp i\in P \, \forall i\}$.
Then $\tdim \Cal Q=\binom n2-1$.
\endproclaim
\demo{Proof}
The dimension of the space of quadrics
containing any $2n+1$ points is
 $\binom n2-1$,
so we need to show that any quadric containing
all but possibly one of the points also contains the last point.
 Recall that  for two
quadratic forms $  Q,  P\in S^2V^*$, that
$\ttrace_{  Q}  P= \Sigma_i   P(v_i,v_i)$ is a well
defined number, where
$\{ v_i\}$ is any $  Q$-orthornormal basis of $V$. Take
$  Q\in S^2V^*$ and let $\hat p_i$ and $\hat z_i$ be
corresponding $  Q$-orthonormal bases. We have
$$
  P(\hat z_0,\hat z_0 ) + \hdots
+  P(\hat z_n,\hat z_n ) =
  P(\hat p_0,\hat p_0 ) + \hdots
+  P(\hat p_n,\hat p_n )
$$
If all the points but perhaps $p_n$ lie on $P$, we see that
$p_n$ must as well.\qed
\enddemo

\bigpagebreak

\noindent{\bf   Relation  to Castelnouvo's lemma 1.12.}
Castelnouvo's lemma says that if $2n+3$ points lie
on an $\{\binom n2 -1\}$-dimensional linear system of quadrics, then in fact
they all lie on a rational normal curve. $\{\binom n2 -1\}$ is the dimension
of the space of quadrics containing a rational normal curve and the point
of Castelnuovo's lemma is   that not only is the space
of the correct dimension, but it is actually a space cutting out
a rational normal curve. Here
we only have $2n+2$ points. (1.11) shows
that the apolarity conditions imply that the $2n+2$ points
always lie on
an $\{\binom n2 -1\}$-dimensional linear system of quadrics.
When one adds the additional hypothesis
on the Cremona images of the points,
the system of quadrics   cuts out a rational normal curve.

\smallpagebreak

\demo{Proof of equivalence of (1.1) and (1.8)}
Write
$$
A
= (p^i_j ) \tag 1.13
$$
and let
$$
\upp j= [ p^0_j\hd p^n_j]. \tag 1.14
$$
Without loss of generality take $\zz i$ as in (1.6) and let
$  Q$ have equation $ \Sigma_i (\xx i)^2=0$.
The Hadamard inverse of $A$ is given by the
coordinates of the $\phi (\upp i)$  up to
the ambiguity of scales which do not effect rank, i.e.
the rank of the Hadamard inverse of $A$ is the dimension of
the span of the $\phi (\upp i)$ plus one. \qed\enddemo

\proclaim{Proposition 1.15, Duality}
 Let $\zz 0\hd\zz n$
and $\upp 0\hd\upp n$ be two sets
 of points
spanning $\pp n$, each set  apolar
with respect to a nondegenerate quadratic form $Q$.
If the
images of the $p_i$ under the Cremona
transform defined by the $z_i$  lie on a $\pp k$, then
images of the $\zz i$ under the Cremona
transform defined by the $\upp i$ also lie on a $\pp k$.
\endproclaim

(1.15) will follow from some remarks on the Gale transform given
in \S 4.

\proclaim{Proposition 1.16}
There exist orthogonal rank $k$ Hadamard inverses for
$m\times m$ complex matrices for all
$k\geq m-\sqrt{\frac{m^2}2+\frac m2 -1}$.
\endproclaim

\demo{Proof}
Let $M=\bcc{m}\otimes \bcc{m}$ and let
$\phi : \ppp M =\pp{m^2-1}\dasharrow \pp{m^2-1}$
be the Cremona transform defined by the
standard coordinates. Write $n=m-1$.

Consider the Segre,
$Seg(\pp n\times\pp n)$,  in $\phi (\ppp M)$
and let $\sigma_k(Seg(\pp n\times\pp n))$ denote its
$k$-th {\it secant variety}, the closure of the union of
all $\pp{k-1}$'s spanned by $k$-ples of points of
the Segre.
Let $Y_k= \phi\inv \sigma_k(Seg(\pp n\times\pp n))$,
so $Y_k$ is the space of matricies with Hadamard inverse
of rank less than or equal to $k$.
Note that $\tdim Y_k = k(2m-k)-1$.

Let $Z\subset\ppp M$ be the space of orthogonal columns,
i.e.
$$
Z= \{ A \, |\, A\upt A \text{ is diagonal} \}.\tag 1.17
$$
$Z$ is a complete intersection of the $\binom m2$ quadrics
$\Sigma_i a^i_ja^i_k=0$ for all $j<k$. ($Z$
is isomorphic to the variety of complete flags.)
$Y_k$ will intersect $Z$ if
$$
k(2m-k)-1 + \binom {m+1}2-1\geq m^2-1, \tag 1.18
$$
i.e., if
$$
k\geq m-\sqrt{\binom {m+1}2 -1}.\qed \tag 1.19
$$
\enddemo

\bigpagebreak

\noindent \S 2. {\it Proof of 1.8.}

Let $A=(p^i_j)$ be as in (1.14) and let $\pp n=\phi (\pp n)$ have
linear coordinates $y^0\hd y^n$ adapted to the image simplex.
Let $H_{ijk}$ ($i<j<k$) denote the hyperplane
defined by the equation
$$
E_{ijk}=
\uup 0i\uup 0 j\uup 0 ky^0+\hdots
 +
\uup ni\uup n j\uup n k y^n. \tag 2.1
$$

Let $q_j=\phi (p_j)$ so
$q^i_j=\uup 0j\hdots \hat p^i_j\hdots \uup nj$,
where the hat denotes omission.
Note that
$$
\tspan \{ q_i,q_j,q_k\}\subseteq H_{ijk}.\tag 2.2
$$
To verify (2.2), by symmetry it is sufficient to verify $q_i\in H_{ijk}$.
$$
\align
E_{ijk}(q_i)&=
\Sigma_l \uup li\uup l j\uup l k(\uup 0i\hdots \hat p^l_i\hdots \uup ni)\\
&=(\uup 0i\hdots   \uup ni)\Sigma_l\uup l j\uup l k\tag 2.3\\
&=(\uup 0i\hdots   \uup ni) \hat Q (\overarrow p_j,\overarrow  p_k) \\
& =0.\endalign
$$
where $\overarrow p_j\in V$ is a  unit vector corresponding to $p_j\in \ppp V$.

Now say
 that the $q_i$   span   a $\pp 2$.
Then
 all the $q_l$'s are in the span of any $q_i,q_j,q_k$
spanning the $\pp 2$. For the moment,
assume we are in the case that there exist two
points, say,
$q_0,q_1$, such that no other $q_i$  lies
on the line between $q_0$ and $q_1$.
(This is
always the case over
$\brrr$, but over $\bccc$ there is the example of the nine flexes on a plane
cubic.)

We claim that the intersection of
$H_{012}\hd H_{01n}$ is at most a
$\pp 1$, which will prove (1.8) in this case.

To see the claim, say there were a linear relation
among $E_{012},E_{013}\hd E_{01n}$,
e.g.
$$
a^2E_{012}+ a^3E_{013}+\hdots + a^n E_{01n}=0. \tag 2.4
$$
The coefficient of $y^j$ in (2.4) is
$$
a^2\uup j0\uup j1\uup j2 + \hdots + a^n\uup j0\uup j1\uup jn =0.
\tag 2.5
$$
This implies (since none of the $\uup ij$ are zero) that
$$
a^2 \uup j2 + \hdots + a^n \uup jn =0
\ \forall j\tag 2.6
$$
i.e. that there is a linear relation among the columns of the matrix
$A$ which is a contradiction.

Now say there are three points, say $q_0,q_1,q_2$ that are colinear
but no other $q_j$ lies on the line they span. We must show that
among the hyperplanes containing the $\pp 2$, that $(n-2)$ of them
are independent. Say not, then for each $3\leq \beta\leq n$ there
must be a relation
$$
\aaa\beta 3E_{013} + \hdots + \aaa\beta nE_{01n} =
a_{\beta}E_{02\beta},
  \tag 2.7
$$
i.e.,
$$
\aaa\beta 3\uup j1\uup j3 + \hdots + \aaa\beta n\uup j1\uup jn =
a_{\beta}\uup j2\uup j\beta
\ \forall j.  \tag 2.8.b
$$
Note that if any of the $a_{\beta}$ are zero we are done by
the above argument.
Similarly there must be a relation
among $H_{023}\hd H_{02n}$ and $H_{013}$ which implies
an equation of the form:
$$
\bb{} 3\uup j2\uup j3 + \hdots + \bb{} n\uup j2\uup jn =
b\uup j1\uup j3. \tag 2.9
$$
Substituting
the right hand side of (2.8.3),...,(2.8.n) into
the left hand side of
(2.9) we   obtain a relation involving
$\uup j1$ in each term which  divides out and we are left with
 a relation
among the rows $p_{\beta}$ and thus a contradiction.

In the event even more points are required, the same argument as
above still works, only one must use more relations.
(To our knowledge there are no known configurations of points
that span a plane with more than three points on each line.)
\qed

\bigpagebreak

\noindent \S 3. {\it Bases and the proof of (1.3)}

\smallpagebreak

\noindent{\bf Definition 3.1}. A {\it base} $\Gamma$ is a set of
$(n+2)$ points in $\bccc\pp n$ in general linear position. The definition is
motivated
by the fact that all such $(n+2)$-ples of points
are projectively equivalent.
Note that there is an $(n-1)$-dimensional linear
system of rational normal curves  through $\Gamma$ which
we will denote $\Cal R_{\Gamma}$ and an
$\{\binom{n+1}2 -2\}$-dimensional linear system of quadric hypersurfaces
through $\Gamma$, which we denote $\Cal Q_{\Gamma}$.

\smallpagebreak

The following is a slight modification  of some
facts in [Con] and [DO]:

\proclaim{Lemma 3.2}  If we fix a
base $\Gamma=\{\zz0\hd\zz n,\upp 0\}$ and
 a hyperplane $H$ with $\Gamma \cap H=\emptyset$,
 then there is a unique quadric
$Q_0$ such that
the $\{\zz i\}$ are   apolar with respect to $Q_0$
and $H$ is the polar hyperplane of $\upp 0$
with respect to $Q_0$.
\endproclaim

\demo{Proof}
Without loss of generality, take $\zz  i$ as vertex points as  in  (1.6) and
take
  $\upp 0 = [1\hd 1]$.
Say $H$ has equation
$\Sigma_i a_i\xx i=0$.  All quadrics for which the $\zz i$ are apolar
are of the form $Q=\Sigma_i\lambda_i (\xx i )^2$ for some constants
$\lambda_i$.
The $p_0$-polar hyperplane of such
a quadric has the equation $\Sigma_i\lambda_i  \xx i$,
so we must have $\lambda_i=a_i$,
uniquely determining $Q_0$.\qed\enddemo

\smallpagebreak

\proclaim{Lemma 3.3}
If we fix a
base $\Gamma=\{\zz0\hd\zz n,\upp 0\}$, a hyperplane
$H$ and $Q_0$ as above, then for all $R \in \Cal R_{\Gamma}$,
the set of $n+1$ points consisting of
$R \cap H$ and $\upp 0$ is apolar with respect to
$Q_0$.
\endproclaim
\demo{Proof}see [DO] Lemma 5, p49.\qed\enddemo

\demo{Proof of (1.3)}
Fixing the $\zz i$ as in (1.6) and $Q_0=\Sigma_i (\xx i)^2$,
i.e. fixing a copy of $SO(n+1)$,
we are free to pick $p_0$ from an open set in $\pp n$,
and then there is a $\pp{n-1}$'s worth of rational normal
curves through $\{ z_i,p_0\}$. We see that the
dimension of complex orthogonal $(n+1)\times (n+1)$ matrices with
rank two Hadamard inverses is $2n-1$.
Finally, taking $Q_0$ as above over $\brrr$ the same
count is still valid. \qed\enddemo

\smallpagebreak

Note that
a rank one Hadamard inverse
is impossible as $\phi$ is one to one off
the hyperplanes that get blown down, and points on blown down
 hyperplanes
correspond to a column vector
with at least one entry equal to zero. (In fact the number of
zeros in a column is the number of such hyperplanes the corresponding
 point lies on.)

\bigpagebreak

\noindent\S 4. {\it Some remarks on the Gale transform}
\smallpagebreak

Everything in this section with the exception
of Version 3 of theorem 1.8 is classical
and explained   in  modern language and greater generality in
[DO], [EP1], and  [EP2]. For our purposes, points will be distinct and in
sufficiently
general linear position so that there is no need to be concerned
with degenerate cases,
and this will enable a simplified exposition.

\smallpagebreak

\noindent{\bf Definition 4.1}. Let $\Gamma$ be a set of
$r+s+2$ points in $\pp r$.  A set $\Gamma '$
of $r+s+2$ points in
$\pp s$ is said to be {\it associated} to $\Gamma$ if when one
chooses coordinates in the respective projective spaces and
  writes the points of $\Gamma$ as the rows of
an $(r+1)\times (r+s+2)$ matrix $A$ and
the points of $\Gamma '$ as the rows of
an $(s+1)\times (r+s+2)$ matrix $B$, that there exists
a   $(r+s+2)\times(r+s+2)$ diagonal matrix $\Lambda$
with nonzero determinant such that $A\Lambda\upt B=0$.
If the point sets are sufficiently nice there is a unique associated
point set. (Of course all this defined up to $PGL(r)$ and $PGL (s)$
actions).

The set $\Gamma'$ is  called the {\it Gale transform} of $\Gamma$.
We will explain association in coordinates.

Let $0\leq i,j\leq r$, and $0\leq \alpha,\beta\leq s$. Write
$\Gamma = \{ \zz i,\upp\alpha \}$, and without loss of generality
(assuming the $\zz i$ are in general linear position) write
$\zz i = [ 0\hd 0,1,0\hd 0]$ where the $1$ is
in the $i$-th position.
Write $\upp\alpha = [\uup 0\alpha\hd\uup r\alpha ]$.
Let $\qq i=[\uup i0\hd\uup is]\in \pp s$ and
let $w_{\alpha}= [ 0\hd 0,1,0\hd 0]\in \pp s$ where the $1$
is in the $\alpha$-th position.

\proclaim{Proposition 4.2} In the situation
above, $\Gamma$ is associated to  $\Gamma '=\{ w_{\alpha},\qq i\}$.
\endproclaim
\demo{Proof}
$$
\pmatrix
1&      & &\uup 00 &\hdots &\uup 0s\\
 &\ddots& &        &\vdots &\vdots\\
 &      &1&\uup r0 & \hdots &\uup rs\endpmatrix
\pmatrix -Id_{r+1}& \\ &Id_{s+1}\endpmatrix
 \pmatrix
\uup 00 &\hdots &\uup 0s\\
        &\vdots &\vdots\\
\uup r0 & \hdots &\uup rs\\
1&      & \\
 &\ddots& \\
  &      &1\endpmatrix
=0 .\qed
$$
\enddemo

\proclaim{Corollary 4.3}   $\tspan\{\upp\alpha \}=\pp k$
if and only if
$\tspan \{\qq i\}=\pp k$.
\endproclaim
\demo{Proof} Row rank equals column rank.\qed\enddemo

\proclaim{Proposition 4.4, Commutativity of
association and Cremona} Let $z^*_i$ be the
image simplex points determined
by the Cremona transform
of $\pp r$in the $\zz i$, which we denote $\phi_z$.
Similarly, let $\qqq *\alpha$
be the simplex points determined
by the Cremona transform of $\pp s$ in the $\qq\alpha$, which we denote
 $\phi_q$.
Then the associated point set of $\Gamma^* =\{ z^*_i,\phi_z(\upp\alpha )\}$
is ${\Gamma '}^* =\{ \phi_q(w_{\alpha}),\qqq *i\}$, i.e.,
${\Gamma '}^* ={\Gamma^*}'$.
\endproclaim
\demo{Proof}
The original set of points $\zz i,\upp\alpha$ yields a matrix
$$
\pmatrix
1&      & &\uup 00 &\hdots &\uup 0s\\
 &\ddots& &        &\vdots &\vdots\\
 &      &1&\uup r0 & \hdots &\uup rs\endpmatrix .
$$
Under $\phi_z$ one gets a point set
$z^*_i,\phi_z (\upp\alpha)$ with matrix
$$
\pmatrix
1&      & &\frac 1{\uup 00} &\hdots &\frac 1{\uup 0s}\\
 &\ddots& &        &\vdots &\vdots\\
 &      &1&\frac 1{\uup r0} & \hdots &\frac 1{\uup rs}\endpmatrix
$$
whose associate is explained above.
But taking transpose commutes with taking Hadamard inverse.\qed
\enddemo

\noindent{\bf Definition 4.5}. A set $\Gamma$ of $2n+2$ points in $\pp n$ is
said to be {\it self associative} if it is associated to itself.

\proclaim{Proposition 4.6} $\Gamma = \{\zz i ,\upp i\}$ is self associative
if and only if the points sets
$\{\zz i\}$ and $\{\upp i\}$ are both
apolar with respect to some quadric $Q$.
\endproclaim
\demo{Proof}
Let  $z_i$ be simplex points as above. We need to show the point
set $\Gamma$
represented by the matrix $(Id, P)$ is equivalent to the point
set represented by the matrix $(Id,\upt P)$, where
$P$ is a matrix whose columns are the entries of the $p_i$,
if and only if the point sets are apolar with respect to a
quadric $Q$.
Let $Q$ also denote  the    $(r+1)\times (r+1)$
matrix representing the quadric.
Let  The point set
$\Gamma'$ consisting of the
columns of $(\upt P, Id)$ is equivalent to
the point set consisting of the
columns of  $(\upt PQ , Q)$. The $p_i$ are $Q$-apolar if
 $\upt PQ =Q P\inv$, in which case
$ (\upt PQ, Q)$ is equivalent to
$(QP\inv, Q)$ which is equivalent to
$(P\inv, Id)$ which is eqivalent to $(Id,P)$.
\qed
\enddemo

In light of (4.6), we can rephrase (1.8) yet again:

\proclaim{Theorem 1.8, Version 3} Let $\Gamma = \{z_i,p_i\}\subset\pp n$
be a self-associated point set. Then $\tspan \{\phi_z(p_i )\}\neq \pp 2$.
  $\tspan \{\phi_z(p_i )\}= \pp 1$ if and only if $\Gamma$ is contained
in a rational normal curve.
\endproclaim

\bigpagebreak

\noindent{\bf Acknowledgements}. It is a pleasure to
thank M. Kontsevich for suggesting the conjecture and
useful comments,
H. Pinkham for many useful conversations and
S. Popescu for    alerting the author to the
relevant literature, interesting remarks, and useful comments.

\bigpagebreak

$$
\text{{\bf References}}
$$

\noindent [Cob1], A.B. Coble, {\it Point sets
and allied Cremona groups I},  Trans. Amer. Math. Soc.
{\bf 16}, (1915), 155-198.

\smallpagebreak

\noindent [Cob2], A.B. Coble, {\it Point sets
and allied Cremona groups II},  Trans. Amer. Math. Soc.
{\bf 17}, (1916), 345-385.

\smallpagebreak

\noindent [Cob3], A.B. Coble, {\it Point sets
and allied Cremona groups III},  Trans. Amer. Math. Soc.
{\bf 18}, (1917), 331-372.

\smallpagebreak

\noindent [Cob4], A.B. Coble, {\it Associated
sets of points},  Trans. Amer. Math. Soc.
{\bf 24}, (1922), 1-20.

\smallpagebreak

\noindent [Con], Conner, {\it Basic systems of rational norm-curves},
American Journal of Mathematics, {\bf 32}(1911), p. 115-176.

\smallpagebreak

\noindent [DO], Dolgachev and Ortland,
{\it Point sets in projective spaces and theta functions},
Ast\'erisque {\bf 165}(1988). 210pp.

\smallpagebreak

\noindent [EP1], D. Eisenbud and S. Popescu,
{\it Gale duality  and
free resolutions of ideals of points.}, preprint.

\smallpagebreak

\noindent [EP2], D. Eisenbud and S. Popescu,
{\it The projective geometry of the Gale transform.}, preprint.

\enddocument